\documentclass[preprint,12pt]{elsarticle}

\usepackage{amssymb}
\usepackage{amsthm}

 \newtheorem{thm}{Theorem}
\newtheorem{lem}[thm]{Lemma}
\newdefinition{rmk}{Remark}
\newproof{pf}{Proof}
\newproof{pot}{Proof of Theorem \ref{principal2}}

\newtheorem{prop}[thm]{Proposition}
\newtheorem{prob}[thm]{Problem}
\newtheorem{coro}[thm]{Corollaire}
\newtheorem{exam}[thm]{Example}

\begin{document}

\begin{frontmatter}

\title{ Skew-symmetric matrices and their principal minors }

\author[authorlabel1]{Abderrahim Boussa\"{\i}ri}
\ead{aboussairi@hotmail.com}
\author[authorlabel1]{Brahim Chergui}
\ead{cherguibrahim@gmail.com}
\address[authorlabel1]{Facult\'e des Sciences A\"{\i}n Chock,
D\'epartement de Math\'ematiques et Informatique, Km 8 route d'El
Jadida, BP 5366 Maarif, Casablanca, Maroc}

\begin{abstract}
Let $V$ be a nonempty finite set and $A=(a_{ij})_{i,j\in V}$ a matrix with
entries in a field $\mathbb{K}$. For a subset $X$ of $V$, we denote by $A[X]$
the submatrix of $A$ having row and column indices in $X$. In this article,
we study the following Problem. Given a positive integer $k$, what is the
relationship between two matrices $A=(a_{ij})_{i,j\in V}$, $%
B=(b_{ij})_{i,j\in V}$ with entries in $\mathbb{K}$ and such that $\det(A%
\left[ X\right] )=\det(B\left[ X\right] )$ for any subset $X$ of $V$ of size
at most $k$ ? The Theorem that we get is an improvement of a result of R.
Lowey \cite{Lw86} for skew-symmetric matrices whose all off-diagonal entries are
nonzero.
\end{abstract}

\end{frontmatter}

\title{Matrices antisym\'{e}triques et leurs mineurs principaux}

% Maintenant la version abrégée en anglais, si présente

\section{Introduction}

Our original motivation comes from the following open problem, called the
Principal Minors Assignement Problem ( PMAP for short) (see \cite{HOZ}).

\begin{prob}
\label{pmap}Find, a necessary and sufficient conditions for a collection of $%
2^{n}$ numbers to arise as the principal minors of a matrix of order $n$ ?
\end{prob}

The PMAP has attracted some attention in recent years. O. Holtz and B.
Strumfels \cite{HOLTZ07} approched this problem algebraically and showed
that a real vector of length $2^{n}$, assuming it strictly satisfies the
Hadamard-Fischer inequalities, is the list of principal minors of some real
symmetric matrix if and only if it satisfies a certain system of polynomial
equations. L. Oeding \cite{Oeding} later proved a more general conjecture of
O. Holtz and B. Strumfels \cite{HOLTZ07}, removing the Hadamard-Fischer
assumption and set-theoretically characterizing the variety of principal
minors of symmetric matrices. Griffin and Tsatsomeros gave an algorithmic
solution to the PMAP \cite{Griffin1,Griffin2}. Their work gives an
algorithm, which, under a certain \textquotedblleft
genericity\textquotedblright\ condition, either outputs a solution matrix or
determines that none exists. Very recently, J. Rising, A. Kulesza, B. Taskar
\cite{Rising} gave an algorithm to solve the PMAP for the symmetric case.
For the general case, the algebraic relations between the principal minors
of a generic matrix of order $n$ are somewhat mysterious. S. Lin and B.
Sturmfels \cite{Lin} proved that the ideal of all polynomial relations among
the principal minors of an arbitrary matrix of order $4$ is minimally
generated by $65$ polynomials of degree $12$. R. Kenyon and R. Pemantle \cite%
{Kenyon} showed that by adding in certain 'almost' principal minors, the
ideal of relations is generated by translations of a single relation.

A natural problem in connection with the PMAP is the following.

\begin{prob}
\label{problem1}What is the relationship between two matrices having equal
corresponding principal minors of all orders ?
\end{prob}

Given two matrices $A$, $B$ of order $n$ with entries in a field $\mathbb{K}$%
, we say that $A$, $B$ are \emph{diagonally similar up to transposition} if
there exist a nonsingular diagonal matrix $D$ such that $B=D^{-1}AD$ or $%
B^{t}=D^{-1}AD$ (where $B^{t}$ is the transpose of $B$). Clearly diagonal
similarity up to transposition preserve all principal minors. D. J. Hartfiel
and R. Loewy \cite{HL}, and then R. Lowey \cite{Lw86} found sufficient
conditions under which diagonal similarity up to transposition is precisely
the relationship that must exist between two matrices having equal
corresponding principal minors. To state the main theorem of \cite{Lw86}, we
need the following notations. Let $V$ be a nonempty finite set and $%
A=(a_{ij})_{i,j\in V}$ be a matrix with entries in a field $\mathbb{K}$ and
having row and column indices in $V$. For two nonempty subsets $X$, $Y$ of $%
V $, we denote by $A\left[ X,Y\right] $ the submatrix of $A$ having row
indices in $X$ and column indices in $Y$. The submatrix $A[X,X]$ is denoted
simply by $A[X]$. The matrix $A$ is\emph{\ irreducible} if for any proper
subset $X$ of $V$, the two matrices $A[X,V\setminus X]$ and $A[V\setminus
X,X]$ are nonzero. The matrix $A$ is\emph{\ HL-indecomposable }(HL after D.
J. Hartfiel and R. Lowey) if for any subset $X$ of $V$ such that $2\leq
|X|\leq |V|-2$, either $A[X,V\setminus X]$ or $A[V\setminus X,X]$ is of rank
at least $2$. Otherwise, it is \emph{HL-decomposable}. The main theorem of
R. Loewy \cite{Lw86} can be stated as follows.

\begin{thm}
\label{loewy} Let $V$ be a nonempty set of size at least $4$ and $%
A=(a_{ij})_{i,j\in V}$, $B=(b_{ij})_{i,j\in V}$ two matrices with entries in
a field $\mathbb{K}$. Assume that $A$ is irreducible and HL-indecomposable.
If $\det (A[X])=\det (B[X])$ for any subset $X$ of $V$, then $A$ and $B$ are
diagonally similar up to transposition.
\end{thm}

For symmetric matrices, the Problem \ref{problem1} has been solved by G. M.
Engel and H. Schneider \cite{Engel}. The following theorem is a directed
consequence of Theorem 3.5 (see \cite{Engel}).

\begin{thm}
\label{symmetric} Let $V$ be a nonempty set of size at least $4$ and $%
A=(a_{ij})_{i,j\in V}$, $B=(b_{ij})_{i,j\in V}$ two complex symmetric
matrices. If $\det (A[X])=\det (B[X])$ for any subset $X$ of $V$ then there
exists a diagonal matrix $D$ with diagonal entries $d_{i}\in \left\{
-1,1\right\} $ such that $B=DAD^{-1}$.
\end{thm}

J. Rising, A. Kulesza, B. Taskar \cite{Rising} gave a simple combinatorial
proof of Theorem \ref{symmetric}. Moreover, they pointed out that the
hypotheses of Theorem \ref{symmetric} can be weakened in sevral special
cases. For symmetric matrices with no zeros off the diagonal, we can improve
this theorem by using the following result. If $A$ is a complex symmetric
matrix with no zeros off the diagonal, then the principal minors of order at
most $3$ of $A$ determine the rest of the principal minors. This result was
obtained by L. Oeding for the generic symmetric matrices (see Remark 7.4,
\cite{Oeding}). Nevertheless, it is not valid for an arbitrary symmetric
matrix. For this, we consider the following example.

\begin{exam}
\label{examplesymm}Let $V:=\{1,...,n\}$ with $n\geq 4$. Consider the
matrices $A_{n}:=(a_{ij})_{i,j\in V}$, $B_{n}:=(b_{ij})_{i,j\in V}$ were $%
a_{i,i+1}=a_{i+1,i}=b_{i,i+1}=b_{i+1,i}=1$ for $i=1,\ldots ,n-1$, $%
a_{n,1}=a_{1,n}=-b_{n,1}=-b_{1,n}=1$ and $a_{ij}=b_{ij}=0$, otherwise. One
can check that $\det (A_{n}\left[ X\right] )=\det (B_{n}\left[ X\right] )$
for any proper subset $X$ of $V$, but $\det (A_{n})\neq \det (B_{n})$.
\end{exam}

Consider now, the skew-symmetric version of Problem \ref{problem1}. We can
ask as for the symmetric matrices if the hypotheses of Theorem \ref{loewy}
can be weakened in special cases. Clearly, two skew-symmetric matrices of
orders $n$ have equal corresponding principal minors of order $3$ if and
only if they differ up to the sign of their off-diagonal entries and they
are not are always diagonally similar up to transposition. Then, it is
necessary to consider principal minors of order $4$. More precisely, we can
suggest the following problem.

\begin{prob}
\label{problem2} Given a positive integer $k\geq 4$, what is the
relationship between two skew-symmetric matrices of orders $n$ having equal
corresponding principal minors of order at most $k$ ?
\end{prob}

Our goal, in this article, is to study this problem for skew-symmetric
matrices whose all off-diagonal entries are nonzero. Such matrices are
called \emph{dense} matrices \cite{Wesp}. A partial answer can be obtained
from the work of G. Wesp \cite{Wesp} about principally unimodular matrices.
A square matrix is\emph{\ principally unimodular } if every principal
submatrix has determinant $0$, $1$ or $-1$. (\cite{BCG}, \cite{Wesp}). G.
Wesp \cite{Wesp} showed that a skew-symmetric dense matrix $%
A=(a_{ij})_{i,j\in V}$ with entries in $\left\{ -1,0,1\right\} $ is
principally unimodular if and only if $\det (A[X])=1$ for any subset $X$ of $%
V$, of size $4$. It follows that if $A$, $B$ are two skew-symmetric dense
matrices have equal corresponding principal minors of order at most $4$,
then they are both principally unimodular or not.

Our main result is the following Theorem which improves Theorem \ref{loewy}
for skew-symmetric dense matrices.

\begin{thm}
\label{principal2} Let $V$ be a nonempty set of size at least $4$ and $%
A=(a_{ij})_{i,j\in V}$, $B=(b_{ij})_{i,j\in V}$ two skew-symmetric matrices
with entries in a field $\mathbb{K}$ of characteristic not equal to $2$.
Assume that $A$ is dense and HL-indecomposable. If $\det (A[X])=\det (B[X])$
for any subset $X$ of $V$ of size at most $4$ then there exists a diagonal
matrix $D$ with diagonal entries $d_{i}\in \left\{ -1,1\right\} $ such that $%
B=DAD$ or $B^{t}=DAD$, in particular, $\det (A)=\det (B)$.
\end{thm}

Note the Theorem \ref{principal2} is not valid for arbitrary
HL-indecomposable skew-symmetric matrices. For this, we consider the
following example.

\begin{exam}
\label{examplskew}Let $V:=\{1,...,n\}$ where $n$ is a even integer such that
$n\geq 6$. Consider the matrices $A_{n}:=(a_{ij})_{i,j\in V}$, $%
B_{n}:=(b_{ij})_{i,j\in V}$ where $%
a_{i,i+1}=b_{i,i+1}=-a_{i+1,i}=-b_{i+1,i}=1$ for $i=1,\ldots ,n-1$, $%
a_{n,1}=-b_{n,1}=-a_{1,n}=b_{1,n}=1$ and $a_{ij}=b_{ij}=0$, otherwise. One
can check that $A_{n}$ and $B_{n}$ are HL-indecomposable and $\det (A_{n}%
\left[ X\right] )=\det (B_{n}\left[ X\right] )$ for any proper subset $X$ of
$V$, but $\det (A_{n})=0$ and $\det (B_{n})=4$.
\end{exam}

Throughout this paper, we consider only matrices whose entries are in a
field $\mathbb{K}$ of characteristic not equal to $2$.

\section{Decomposability of Matrices}

Let $A=(a_{ij})_{i,j\in V}$ \ be a matrix. A subset $X$ of $V$ is a \emph{%
HL-clan} of $A$ if both of matrices $A\left[ X,V\setminus X\right] $ and $A%
\left[ V\setminus X,X\right] $ have rank at most $1$. By definition, the
complement of an HL-clan is an HL-clan. Moreover, $\emptyset $, $V$,
singletons $\{x\}$ and $V\setminus \{x\}$ (where $x\in V$) are HL-clans of $A
$ called \emph{trivial }HL-clans. It follows that a matrix is
HL-indecomposable if all its HL-clans are trivial.

From clan's definition for 2-structures \cite{EHR ROZ}, we can introduce the
concept of clan for matrices which is stronger that of HL-clan. A subset $I$
of $V$ is a \emph{clan} of $A$ if for every $i,j\in I$ and $x\in V\setminus
I $, $a_{xi}=a_{xj}$ and $a_{ix}=a_{jx}$. For example, $\emptyset$, $\left\{
x\right\} $ where $x\in V$ and $V$ are clans of $A$ called \emph{trivial }%
clans. A matrix is \emph{indecomposable} if all its clans are trivial.
Otherwise, it is \emph{decomposable}.

In the next proposition, we present some basic properties of clans that can
be deduced easily from the definition.

\begin{prop}
\label{intervalpropriete} Let $A=(a_{ij})_{i,j\in V}$ be a matrix and let $%
X,Y$ be two subsets of $V$.
\begin{item}
\item[i)] If $X$ is a clan of $A$ and $Y$ is a clan of $A[X]$ then $Y$ is a clan of $A$.
\item[ii)] If $X$ is a clan of $A$ then $X\cap Y$ is a clan of $A\left[  Y\right]  $.
\item[iii)]If $X$ and $Y$ are clans of $A$ then $X\cap Y$ is a clan of $A$.
\item[iv)] If $X$ and $Y$ are clans of $A$ such that $X\cap Y\neq\emptyset$, then $X\cup
Y$ is a clan of $A$.
\item[v)] If $X$ and $Y$ are clans of $A$ such that $X\setminus Y\neq\emptyset$, then
$Y\setminus X$ is a clan of $A$.

\end{item}

\end{prop}

Let $A=(a_{ij})_{i,j\in V}$ be a matrix and let $X,Y$ be two nonempty
disjoint subsets of $V$. If for some $\alpha\in\mathbb{K}$, we have $%
a_{xy}=\alpha$ for every $x\in X$ and every $y\in Y$ then we write $%
A_{(X,Y)}=\alpha$. Clearly, if $X,Y$ are two nonempty disjoint clans of $A$
then there is $\alpha \in\mathbb{K}$ such that $A_{(I,J)}=\alpha$.

A \emph{clan-partition} of $A$ is a partition $\mathcal{P}$ of $V$ such that
$X$ is a \emph{clan} of $A$ for every $X\in\mathcal{P}$.

The following lemma gives a relationship between decomposability and
HL-decomposability.

\begin{lem}
\label{indecoet normalisation}

Let $V$ be a finite set of size at least $4$ and $A=(a_{ij})_{i,j\in V}$ a
matrix. We assume that there is $v\in V$ and $\lambda$, $\kappa$ $\in\mathbb{%
K}\setminus\{0\}$ such that $a_{vj}=\lambda$ et $a_{jv}=\kappa$ for $j\neq v$%
. Then $A$ is HL-decomposable if and only if $A[V\setminus\{v\}]$ is
decomposable.
\end{lem}

\begin{pf}
Clearly $V\setminus\{v\}$ is a clan of $A$, it follows by Proposition \ref%
{intervalpropriete} that a clan of $A[V\setminus\{v\}]$ is a clan of $A$ and
hence it is a HL-clan of $A$. Therefore, if $A[V\setminus\{v\}]$ is
decomposable then $A$ is HL-decomposable. Conversely, suppose that $A$ is
HL-decomposable and $I$ is a nontrivial HL-clan of $A$. As $V\setminus I$ is
also a nontrivial HL-clan of $A$ then, by interchanging $I\ $and $V\setminus
I$, we can assume that $I\subseteq V\setminus\left\{ v\right\} $. We will
show that $I$ is a clan of $A[V\setminus\{v\}]$. For this, let $i,j\in I$
and $k\in(V\setminus\left\{ v\right\} )\setminus I$. The matrix $%
A[\{v,k\},\{i,j\}]$ (resp. $A[\{i,j\},\{v,k\}]$) is a submatrix of $%
A[V\setminus I,I]$ (resp. $A[I,V\setminus I]$). As, $I$ and $V\setminus I$
are HL-clans of $A$ then both of matrices $A[V\setminus I,I]$ and $%
A[I,V\setminus I]$ have rank at most $1$. It follows that $\det(A\left[
\left\{ v,k\right\} ,\left\{ i,j\right\} \right] )=0$, $\det(A\left[ \left\{
i,j\right\} ,\left\{ v,k\right\} \right] )=0$ and so, $a_{ki}=a_{kj}$ and $%
a_{ik}=a_{jk}$ . We conclude that $I$ is a nontrivial clan of $%
A[V\setminus\{v\}]$.
\end{pf}

\begin{lem}
\label{memeclan} Let $V$ be a finite set, $A=(a_{ij})_{i,j\in V}$ a matrix
and $D=(d_{ij})_{i,j\in V}$ a nonsingular diagonal matrix. The matrices $A$
and $DAD$ have the same HL-clans. In particular, $A$ is HL-indecomposable if
and only if $DAD$ is HL-indecomposable.
\end{lem}

\begin{pf}
Let $X$ be a subset of $V$. We have the following equalities :

$(DAD)\left[ V\setminus X,X\right] =(D\left[ V\setminus X\right] )(A\left[
V\setminus X,X\right] )(D\left[ X\right] )$

$(DAD)\left[ X,V\setminus X\right] =(D\left[ X\right] )(A\left[ X,V\setminus
X\right] )(D\left[ V\setminus X\right] )$

But, the matrices $D\left[ X\right] $ and $D\left[ V\setminus X\right] $ are
nonsingular, then $(DAD)\left[ V\setminus X,X\right] $ and $A\left[
V\setminus X,X\right] $ (resp. $(DAD)\left[ X,V\setminus X\right] $ and $(A%
\left[ X,V\setminus X\right] )$ have the same rank. Therfore, the matrices $%
A $ and $DAD$ have the same HL-clans.
\end{pf}

\section{Separable matrices}

We introduce here the notion of separable matrix that comes from
clan-decomposability of 2-structures (see [3] ). A skew-symmetric matrix $%
A=(a_{ij})_{i,j\in V}$ is \emph{separable} if it has a clan whose complement
is also a clan. Otherwise, it is \emph{inseparable}.

\begin{lem}
\label{lemnouv} Let $A=(a_{ij})_{i,j\in V}$ be a separable matrix and $%
\left\{ X,Y\right\} $ a clan-partition of $A$. If there is $x\in X$ such
that $A\left[ V\setminus\left\{ x\right\} \right] $ is inseparable, then $%
X=\left\{ x\right\} $ and $Y=V\setminus\left\{ x\right\} $.
\end{lem}

\begin{pf}
We assume, for contradiction, that $X\neq\left\{ x\right\} $. From
Proposition \ref{intervalpropriete}, $\{X\cap ( V\setminus\left\{
x\right\}), Y\} $ is a clan-partition of $A\left[ V\setminus\left\{
x\right\} \right] $ but this contradicts the fact that $A\left[
V\setminus\left\{ x\right\} \right] $ is an inseparable matrix. It follows
that $X=\left\{ x\right\} $ and $Y=V\setminus\left\{ x\right\} $.
\end{pf}

\begin{lem}
\label{decominseparable}

Let $V$ be a finite set with $\left\vert V\right\vert \geq5$ , $%
A=(a_{ij})_{i,j\in V}$ a skew-symmetric matrix and $Y$ a subset of $V$ such
that $\left\vert Y\right\vert \geq2$. If $C$ is a clan of $A$ such that $%
|C\cap Y|=1$ and $V=C\cup Y$, then $A$ is inseparable if and only if $A[Y]$
is inseparable.

\end{lem}

\begin{pf}

Let $C\cap Y:=\{y\}$. Suppose that $A[Y]$ is separable and let $%
\{Z,Z^{\prime }\}$ a clan-partition of $A[Y]$. By interchanging $Z$ by $%
Z^{\prime}$, we can assume that $y\in Z$. It is easy to see that $\left\{
Z\cup C,Z^{\prime}\right\} $ is a clan-partition of $A$ and hence $A$ is
separable. Conversely, suppose that $A[Y]$ is inseparable. Let $W$ be a
subset be of $V$ containing $Y$ and such that $A[W]$ is inseparable. Choose $%
W$ so that its cardinality is maximal. To prove that $A$ is inseparable, we
will show that $W=V$. For contradiction, suppose that $W\neq V$ and let $%
v\in V\setminus W$. By maximality of $\left\vert W\right\vert $, the matrix $%
A[W\cup\{v\}]$ is separable and by Lemma \ref{lemnouv}, $\left\{ \left\{
v\right\} ,W\right\} $ is a clan-partition of $A[W\cup\{v\}]$. Furthermore,
by Proposition \ref{intervalpropriete}, $C\cap(W\cup\{v\})$ is a clan of $%
A[W\cup\{v\}]$. Moreover, as $v\in C$ because $C\cup Y=V$, $Y\subseteq W$
and $v\in V\setminus W$, it follows that $v\in(C\cap(W\cup\{v\}))\setminus W$
and hence $(C\cap(W\cup\{v\}))\setminus W\neq\emptyset$. By applying
Proposition \ref{intervalpropriete}, the subset $W\setminus(C\cap\left(
W\cup\left\{ v\right\} \right) )=W\setminus C$ is a clan of $A\left[
W\cup\left\{ v\right\} \right] $. Thus $\left\{ W\setminus C,C\cap W\right\}
$ is a clan-partition of $A[W]$ but this contradicts the fact that $A[W]$ is
inseparable.

\end{pf}

As a consequence of this result, we obtain the following corollary.

\begin{coro}
\label{moon} Let $V$ be a finite set of size $n$ with $n\geq5$ and $%
A=(a_{ij})_{i,j\in V}$ an inseparable skew-symmetric matrix. Then, there is $%
x\in V$ such that the matrix $A[V\setminus\{x\}]$ is inseparable.
\end{coro}

\begin{pf}
First, suppose that the matrix $A$ is decomposable. Let $I$ be a nontrivial
clan of $A$ and $y\in I$. Put $Y:=(V\setminus I)\cup\{y\}$. Clearly, $Y\cap
I=\left\{ y\right\} $ and $V=I\cup Y$. By Lemma \ref{decominseparable}, the
matrix $A[Y]$ is inseparable. Let $x\in I\setminus\{y\}$. \ As $I$ is a clan
of $A$, then by Proposition \ref{intervalpropriete}, $I\setminus\{x\}$ is a
clan of $A[V\setminus\{x\}]$. In addition, $Y\cap(I\setminus\{x\})=\{y\}$
and $(I\setminus\{x\})\cup Y=V\setminus\{x\}$, then, by applying again Lemma %
\ref{decominseparable}, we deduce that the matrix $A\left[ V\setminus\left\{
x\right\} \right] $ is necessarily inseparable.

Suppose, now, that the matrix $A$ is indecomposable.

Firstly, by proceeding as in the proof of Theorem 6.1 (pp. 93 \cite{EHR ROZ}%
), we will show that there is a subset $I$ of $V$ such that $3\leq|I|\leq
n-1 $ and $A[I]$ is inseparable. Assume the contrary and let $x\in V$ . The
matrix $A[V\setminus\{x\}]$ is, then, separable. Let $\{X,X^{\prime}\}$ be a
clan-partition of $A[V\setminus\{x\}]$. Without loss of generality, we can
assume that $|X^{\prime}|\geq|X|$. It follows that $|X^{\prime}|\geq2$
because $|V|\geq5$. As $3\leq|X^{\prime}\cup\left\{ x\right\} |\leq n-1$,
the submatrix $A[X^{\prime}\cup\left\{ x\right\} ]$ is separable. Let $%
\{Y,Y^{\prime}\}$ be a clan-partition of $A[X^{\prime}\cup\{x\}]$ such that $%
x\in Y$. Let $A_{(Y^{\prime},Y)}=a$ and $A_{(X^{\prime},X)}=b$. Since $%
Y^{\prime}\subseteq X^{\prime}$, we have $A_{(Y^{\prime},X)}=b$. In
addition, $X\cup Y=V\setminus$ $Y^{\prime}$, then $Y^{\prime}$ is a clan of $%
A$. It follows that $Y$ is a singleton because $A$ is indecomposable. Let $%
Y^{\prime}:=\{y\}$. If $a=b$ then $A_{(y,X)}=b=a$ and $A_{(y,Y)}=a$. It
follows that $A_{(y,V\setminus\left\{ y\right\} )}=a$ and this contradicts
the fact that $A$ is indecomposable. So $a\not =b$. If$\ A_{(x,X)}=b$, then $%
A_{(V\setminus X,X)}=b$ because $A_{(X^{\prime},X)}=b$ and $V\setminus
X=X^{\prime}\cup\left\{ x\right\} $ and this contradicts that $A$ is
indecomposable. It follows that there is $z\in X$ such that $a_{xz}\not =b$.
Now, by hypothesis, $A[x,y,z]$ is separable, so $a_{zx}=a$. Furthermore, $%
(X^{\prime}\cup\left\{ x\right\} )\setminus \{x,y\}\neq\emptyset$ because $%
|X^{\prime}\cup\left\{ x\right\} |\geq3$. Let $v\in(X^{\prime}\cup\left\{
x\right\} )\setminus\left\{ x,y\right\} =X^{\prime}\setminus\left\{
y\right\} $ . We can choose $v$ such that $a_{vx}\not =a$, because otherwise
$X$ would be a nontrivial clan of $A$ and this contradicts the fact that $A$
is indecomposable. The submatrix $A[x,z,v]$ is also separable, so $a_{vx}=b$%
. It follows that the submatrix $A\left[ x,y,z,v\right] $ is indecomposable,
and then it is inseparable, but this contradicts our assumption. Thus there
is a subset $I$ of $V$ such that $3\leq\left\vert I\right\vert \leq n-1$ and
$A\left[ I\right] $ is inseparable.

To continue the proof, we choose a subset $I$ of $V$ with maximal
cardinality such that $3\leq\left\vert I\right\vert \leq n-1$ and $A\left[ I%
\right] $ is inseparable. We will show that $\left\vert I\right\vert =n-1$.
We assume, for contradiction, that $\left\vert I\right\vert <n-1$. Let $x\in
V\setminus I$. By maximality of $\left\vert I\right\vert $, the matrix $%
A[I\cup\{x\}]$ is separable. Let $\{Z,W\}$ a clan-partition of $A[I\cup
\{x\}]$ such that $x\in Z$. As $A\left[ I\right] $ is inseparable, then by
Lemma \ref{lemnouv}, we have $Z=\{x\}$ and $W=I$. It follows that $I$ is a
proper clan of $A$ which contradicts the fact that $A$ is indecomposable.

\end{pf}

\section{HL-equivalent matrices}

Let $V$ be a finite set and $k$ a positive integer. Two matrices $%
A=(a_{ij})_{i,j\in V}$, $B=(b_{ij})_{i,j\in V}$ are $(\leq k)-$\emph{%
HL-equivalents} if for any subset $X$ of $V$, of size at most $k$, we have $%
\det(A\left[ X\right] )=\det(B\left[ X\right] )$.

In the rest of this paper, we will consider only skew symmetric dense
matrices.

It is clear that if two skew symmetric dense matrices $A=(a_{ij})_{i,j\in V}$%
, $B=(b_{ij})_{i,j\in V}$ are $(\leq2)-$HL-equivalents then, $b_{xy}=a_{xy}$
or $b_{xy}=-a_{xy}$ for all $x\neq y\in V$. We define on $V$ two equivalence
relations $\mathcal{E}_{A,B}$ and $\mathcal{D}_{A,B}$ as follows. For $x\neq
y\in V$, $x\mathcal{E}_{A,B}y$ (resp. $x\mathcal{D}_{A,B}y$) if there is a
sequence $x_{0}=x,x_{1},\ldots,x_{n}=y$ such that $%
a_{x_{i}x_{i+1}}=b_{x_{i}x_{i+1}}$ (resp. $%
a_{x_{i}x_{i+1}}=-b_{x_{i}x_{i+1}} $) for $i=0,\ldots,n-1$. Let $x\in V$, we
denote by $\mathcal{E}_{A,B} \langle x \rangle$ (resp. $\mathcal{D}_{A,B}
\langle x \rangle$) the equivalence class of $x$ for $\mathcal{E}_{A,B}$ (
resp. for $\mathcal{D}_{A,B}$)

Let $V$ be a finite set and $A=(a_{ij})_{i,j\in V}$ a dense skew-symmetric
matrix. Let $V^{\infty}:=V\cup\left\{ \infty\right\} $ and $A^{\infty
}:=(a_{ij})_{i,j\in V^{\infty}}$ were $a_{\infty\infty}=0$, $a_{\infty j}=1$
and $a_{j\infty}=-1$ for $j\in V$.

\begin{lem}
\label{ihemorphe 3sommets}

Let $V=\left\{ i,j,k\right\} $ and $A=(a_{ij})_{i,j\in V}$, $%
B=(b_{ij})_{i,j\in V}$ two dense skew-symmetric matrices such that $%
a_{ij}=b_{ij}$, $a_{ki}=b_{ik}$ and $a_{kj}=b_{jk}$. If $\det\left(
A^{\infty}\right) =\det\left( B^{\infty}\right) $ then $a_{ik}=a_{jk}$ and $%
b_{ik}=b_{jk}$.

\end{lem}

\begin{pf}

We have the following equalities :

$\det(A^{%
\infty})=a_{ij}^{2}-2a_{ij}a_{ik}+2a_{ij}a_{jk}+a_{ik}^{2}-2a_{ik}a_{jk}+a_{jk}^{2}
$,

$\det(B^{%
\infty})=b_{ij}^{2}-2b_{ij}b_{ik}+2b_{ij}b_{jk}+b_{ik}^{2}-2b_{ik}b_{jk}+b_{jk}^{2}
$.

But, by hypothesis, $b_{ij}=a_{ij}$, $b_{ik}=a_{ki}=-a_{ik}$ and $%
b_{jk}=a_{kj}=-a_{jk}$, then we have,

$\det(B^{\infty})=a_{ij}^{2}+2a_{ij}a_{ik}-2a_{ij}a_{jk}+a_{ik}
^{2}-2a_{ik}a_{jk}+a_{jk}^{2}$.

It follows that, if $\det(A^{\infty})=\det(B^{\infty})$ then $a_{jk}=a_{ik}$
because $a_{ij}\neq0$.

Moreover, $b_{ik}=-a_{ik}=-a_{jk}=-b_{kj}=b_{jk}$. Then, we have also $%
b_{ik}=b_{jk}$.

\end{pf}

\begin{lem}
\label{lopez} Let $V$ be a finite set with $|V|\geq3$ and $%
A=(a_{ij})_{i,j\in V}$, $B=(b_{ij})_{i,j\in V}$ two dense skew-symmetric
matrices. If $A^{\infty}$and $B^{\infty}$ are $(\leq4)$ -HL-equivalent then
the equivalence classes of $\mathcal{E}_{A,B}$ and those of $\mathcal{D}%
_{A,B}$ are clans of $A$ and of $B$.

\end{lem}

\begin{pf}

First, note that $\mathcal{D}_{A,B}=\mathcal{E}_{A,B^{t}}$ and the matrices $%
B^{\infty }$, $(B^{t})^{\infty }$ are $(\leq 4)$ -HL -equivalents. So, by
interchanging $B$ and $B^{t}$, it suffices to show that an equivalence class
$C$ of $\mathcal{E}_{A,B}$ is a clan of $A$ and of $B$. For this, let $x\neq
y\in C$ and $z\in V\setminus C$. By definition of $\mathcal{E}_{A,B}$, there
exists a sequence $x_{0}=x,x_{1},\ldots ,x_{r}=y$ such that $%
a_{x_{i}x_{i+1}}=b_{x_{i}x_{i+1}}$for $i=0,\ldots ,r-1$. As $z\notin C$,
then for $j=0,\ldots ,r$, we have $a_{zx_{j}}=$ $-b_{zx_{j}}=b_{x_{j}z}$. It
follows, from Lemma \ref{ihemorphe 3sommets}, that $a_{zx_{i}}=a_{zx_{i+1}}$%
and $b_{zx_{i}}=b_{zx_{i+1}}$ for $i=0,\ldots ,r-1$. Therefore $a_{zx}=a_{zy}
$, $b_{zx}=b_{zy}$ and hence $C$ is a clan of $A$ and of $B$.

\end{pf}

\begin{prop}
\label{plusieursclasse}

Let $V$ be a finite set with $|V|\geq3$ and $A=(a_{ij})_{i,j\in V}$, $%
B=(b_{ij})_{i,j\in V}$ two dense skew-symmetric matrices such that $%
A^{\infty}$ and $B^{\infty}$ are $(\leq4)$ -HL-equivalents. If $A$ is
inseparable, then, $\mathcal{D}_{A,B}$ or $\mathcal{E}_{A,B}$ has at least
two equivalence classes.

\end{prop}

\begin{pf}

The result is trivial for $|V|=2$ or $|V|=3$.

For $|V|=4$, let $V:=\{i,j,k,l\}$. Suppose that $V$ is an equivalence class
for $\mathcal{D}_{A,B}$ and for $\mathcal{E}_{A,B}$. We will show that $A$
is separable. Without loss of generality, we can assume that $a_{il}=-b_{il}$%
, $a_{lk}=-b_{lk}$, $a_{kj}=-b_{kj}$, and $a_{ki}=b_{ki}$, $a_{ij}=b_{ij}$, $%
a_{jl}=b_{jl}$. As $A^{\infty}$ and $B^{\infty}$ are $(\leq4)$
-HL-equivalents, then $\det(A^{\infty}\left[ \left\{ \infty\right\} \cup X%
\right] )=\det(B^{\infty}\left[ \left\{ \infty\right\} \cup X\right] )$ for $%
X=\{i,j,k\}$, $X=\{i,j,l\}$, $X=\left\{ i,k,l\right\} $ and $X=\left\{
i,k,l\right\} $. It follows that $4a_{ij}a_{jk}-4a_{ik}a_{jk}=0$, $%
-4a_{ij}a_{il}-4a_{il}a_{jl}=0$, $4a_{ik}a_{kl}-4a_{ik}a_{il}=0$, $%
-4a_{jk}a_{jl}-4a_{jl}a_{kl}=0$ and then, $a_{ij}=a_{ik}$, $a_{ij}=-a_{jl}$,
$a_{kl}=a_{il}$, $a_{jk}=-a_{kl}$. Furthermore, $\det(A)=\det(B)$, so $%
4a_{kl}a_{ij}\left( a_{ik}a_{jl}-a_{il}a_{jk}\right) =0$ and consequently $%
a_{ik}a_{jl}-a_{il}a_{jk}=0$. We conclude that $a_{ij}^{2}-a_{kl}^{2}=0$. If
$a_{ij}=a_{kl}$ then $a_{ij}=a_{ik}=a_{il}=a_{kl}=a_{kj}=a_{lj}$ and hence $%
\{i,j\}$ and $\{k,l\}$ are two clans of $A$. If $a_{ij}=-a_{kl}$ then $%
a_{ij}=a_{ik}=a_{li}=a_{lj}=a_{jk}=a_{lk}$. So $\{j,k\}$ and $\{i,l\}$ are
two clans of $A$ and hence $A$ is separable.

As in the proof of Proposition 3 (see \cite{BILT}), we will continue by
induction on $|V|$ for $|V|\geq 5$. Suppose, by contradiction, that $A$ is
inseparable and that $V$ is an equivalence class for $\mathcal{D}_{A,B}$ and
for $\mathcal{E}_{A,B}$. By Corollary \ref{moon}, there is $v\in V$ such
that the submatrix $A[V\setminus \{v\}]$ is inseparable. Let $\mathcal{D}%
^{\prime }:=\mathcal{D}_{A\left[ V\setminus \left\{ v\right\} \right] ,B%
\left[ V\setminus \left\{ v\right\} \right] }$ and $\mathcal{E}^{\prime }:=%
\mathcal{E}_{A\left[ V\setminus \left\{ v\right\} \right] ,B\left[
V\setminus \left\{ v\right\} \right] }$. By induction hypothesis and without
loss of generality, we can assume that $\mathcal{D}^{\prime }$ has at least
two equivalence classes. As $\mathcal{E}\langle v\rangle =V$, then there
exist $x\neq v$ such that $a_{xv}=b_{xv}$. Moreover, the matrix $%
A[V\setminus \{v\}]$ is inseparable, then there exist $y_{1}\in \mathcal{D}%
^{\prime }\langle x\rangle $ and $y_{2}\in (V\setminus \left\{ v\right\}
)\setminus \mathcal{D}^{\prime }\langle x\rangle $ such that $a_{xv}\neq $ $%
b_{y_{1}y_{2}}$. Since $\mathcal{D}^{\prime }\langle y_{2}\rangle
\varsubsetneq V$ and $\mathcal{D}_{A,B}\langle y_{2}\rangle =V$, there
exists $z_{2}\in \mathcal{D}^{\prime }\langle y_{2}\rangle $ and $w_{1}\in
V\setminus \mathcal{D}^{\prime }\langle y_{2}\rangle $ such that $%
a_{z_{2}w_{1}}\neq b_{z_{2}w_{1}}$. But, $\mathcal{D}^{\prime }\langle
y_{2}\rangle =\mathcal{D}^{\prime }\langle z_{2}\rangle $, so for any $w\in
(V\setminus \left\{ v\right\} )\setminus \mathcal{D}^{\prime }\langle
y_{2}\rangle $ we have $a_{z_{2}w}=b_{z_{2}w}$. We have necessarily, $w_{1}=v
$ and hence $a_{z_{2}v}\neq b_{z_{2}v}$. Furthermore, $\mathcal{D}^{\prime
}\langle y_{2}\rangle =\mathcal{D}^{\prime }\langle z_{2}\rangle $ and $%
\mathcal{D}^{\prime }\langle y_{2}\rangle \cap \mathcal{D}^{\prime }\langle
x\rangle =\emptyset $ because $y_{2}\in (V\setminus \left\{ v\right\}
)\setminus \mathcal{D}^{\prime }\langle x\rangle $, then $z_{2}$ $\notin
\mathcal{D}^{\prime }\langle x\rangle $ and hence $a_{z_{2}x}=b_{z_{2}x}$.
To obtain the contradiction, we will apply Lemma \ref{ihemorphe 3sommets} to
$X=\left\{ v,x,z_{2}\right\} $. By the foregoing, we have $a_{z_{2}v}\neq
b_{z_{2}v}$, $a_{xv}=b_{xv}$ and $a_{z_{2}x}=b_{z_{2}x}$. But, by Lemma \ref%
{lopez}, $\mathcal{D}^{\prime }\langle x\rangle $ and $\mathcal{D}^{\prime
}\langle y_{2}\rangle $ are disjoint clans of $A[V\setminus \{v\}]$, so $%
a_{y_{1}y_{2}}=a_{xz_{2}}$ because $x,y_{1}\in \mathcal{D}^{\prime }\langle
x\rangle $ and $y_{2},z_{2}\in \mathcal{D}^{\prime }\langle y_{2}\rangle $.
It follows that $a_{xz_{2}}\neq a_{xv}$ because $a_{xv}\neq a_{y_{1}y_{2}}$
and this contradicts Lemma \ref{ihemorphe 3sommets}.

\end{pf}

\section{Proof of the main Theorem}

By using the proposition \ref{plusieursclasse} and Lemma \ref{lopez}, we
obtain the following result which is the key of the proof of the main
Theorem.

\begin{prop}
\label{bilt} Let $V$ be a finite set with $|V|\geq2$ and $A=(a_{ij})_{i,j\in
V}$, $B=(b_{ij})_{i,j\in V}$ two dense skew-symmetric matrices such that $%
A^{\infty}$ and $B^{\infty}$ are $(\leq4)$- HL-equivalent. If $A$ is
indecomposable then $B=A$ or $B=A^{t}$.

\end{prop}

\begin{pf}

Assume that the matrix $A$ is indecomposable, then it is inseparable. By
using proposition \ref{plusieursclasse} and by interchanging $A$ and $A^{t}$%
, we can assume that $\mathcal{D}_{A,B}$ at least two equivalence classes.
By Lemma \ref{lopez} these classes are clans of $A$. So these are singletons
because $A$ is indecomposable. Consequently $B=A$.

\end{pf}

\begin{pot}
Let $u\in V$ and $D:=(d_{ij})_{i,j\in V}$ (resp. $D^{\prime}:=(d_{ij}^{%
\prime })_{i,j\in V}$) the diagonal matrix such that $d_{uu}=1$ and $d_{zz}=%
\frac {1}{a_{uz}}$ (resp. $d_{uu}^{\prime}=1$ and $d_{zz}^{\prime}=\frac{1}{%
b_{uz}}$) if $z\neq u$. We consider the two matrices $\widehat{A}:=DAD:=(%
\widehat{a}_{ij})_{i,j\in V}$ and $\widehat{B}:=D^{\prime}BD^{\prime}:=(%
\widehat{b}_{ij})_{i,j\in V}$. As the matrix $A$ is HL-indecomposable, then
by lemma \ref{memeclan}, the matrix $\widehat{A}$ is also HL-indecomposable.
However $\widehat{a}_{vj}=\widehat{b}_{vj}=1$ for $j\neq v$, then by Lemma %
\ref{indecoet normalisation}, $\widehat{A}\left[ V\setminus\left\{ v\right\} %
\right] $ is indecomposable. We will apply Proposition \ref{bilt}. For this,
let $X$ be a subset of $V$ such that $\left\vert X\right\vert \leq4$. We
have $\det(\widehat {A}\left[ X\right] )=\left( \det D\left[ X\right]
\right) ^{2}\det(A\left[ X\right] )$. Similarly, $\det(\widehat{B}\ \left[ X%
\right] )=\left( \det D^{\prime}\left[ X\right] \right) ^{2}\det(B\left[ X%
\right] )$. But for all $z\neq u$, $b_{uz}=a_{uz}$ or $b_{uz}=-a_{uz}$, then
$\left( \det D\left[ X\right] \right) ^{2}=\left( \det D^{\prime}\left[ X%
\right] \right) ^{2}$. Moreover, by hypothesis, $\det(A[X])=\det(B[X])$,
then $\det(\widehat{A}\left[ X\right] )=\det(\widehat{B}\ \left[ X\right] )$
and hence the matrices $\widehat{A}$, $\widehat{B}\ $ are $(\leq4)$%
-HL-equivalents. Now, by Proposition \ref{bilt}, $\widehat{A}\left[
V\setminus\left\{ v\right\} \right] =\widehat{B}\left[ V\setminus\left\{
v\right\} \right] $ or $\widehat{A}\left[ V\setminus\left\{ v\right\} \right]
=(\widehat{B}\left[ V\setminus\left\{ v\right\} \right] )^{t}$. If $\widehat{%
A}\left[ V\setminus\left\{ v\right\} \right] =\widehat {B}\left[
V\setminus\left\{ v\right\} \right] $ then $\widehat{A}=\widehat{B}$ \ and
hence $B=(D^{\prime})^{-1}DA(D^{\prime})^{-1}D$. Suppose that $\widehat{A}%
\left[ V\setminus\left\{ v\right\} \right] =(\widehat {B}\left[
V\setminus\left\{ v\right\} \right] )^{t}$ and let $\Delta:=(%
\delta_{ij})_{i,j\in V}$ the diagonal matrix such that $\delta _{uu}=1$ and $%
\delta_{zz}=-1$ if $z\neq u$. Clearly, $\Delta\widehat{A}\Delta=(\widehat{B}%
)^{t}$. It follows that $B^{t}=(D^{\prime})^{-1}D\Delta
A(D^{\prime})^{-1}D\Delta$. However for all $z\in V$, $d_{zz}=d_{zz}^{%
\prime} $ or $d_{zz}=-d_{zz}^{\prime}$. then $(D^{\prime})^{-1}D$ and $%
(D^{\prime})^{-1}D\Delta$ are diagonal matrices with diagonal entries in $%
\left\{ -1,1\right\} $.

\end{pot}

\section{Some remarks about Theorem \protect\ref{principal2}}

Let $A$ be a skew-symmetric matrix. Clearly, the principal minors of order $%
2 $ determine the off diagonal entries of $A$ up to sign. So, generically
there are at most finitely many skew-symmetric matrices with equal
corresponding principal minors as a given matrix, and one should expect
Theorem \ref{principal2} to hold for sufficiently generic skew-symmetric
matrices. Nevertheless, for dense skew-symmetric matrix, Theorem \ref%
{principal2} is the best one. To see this, let $A=(a_{ij})_{i,j\in V}$ be an
HL-decomposable matrix and $X$ a nontrivial HL-clan of $A$. Consider the
matrix $B=(b_{ij})_{i,j\in V}$ such that $b_{ij}=-a_{ij}$ if $i$,$j\in X$
and $b_{ij}=a_{ij}$ otherwise.\ By adapting the proof of Lemma 5 (see \cite%
{HL} ), we can prove that $A$ and $B$ have the same principal minors, but
they are not diagonally similar up to transposition.

We can ask if Theorem \ref{principal2} can be obtained from Theorem \ref%
{loewy} via specialization to skew-symmetric matrices. For this, we must
prove that for a dense skew-symmetric matrix, the principal minors of order
at most $4$ determine the rest of its principal minors. One way to do this
is to give an expression of the determinant of a skew-symmetric matrix of
order at least $5$ from its principal minors. More precisely, let $M$ be a
generic skew-symmetric matrix of order $n$ where $n\geq 5$ and the entries $%
x_{ij}$ with $i$ $<j$ are indeterminates. Let $R=\mathbb{K}[x_{ij}|1\leq
i<j\leq n]$ be the polynomial ring generated by these indeteminates. The
problem is to find an expression of $\det \left( M\right) $ from the
collection of the principal minors $\det \left( M\left[ \alpha \right]
\right) $ where $\alpha $ is a subset of $\left\{ 1,\ldots ,n\right\} $ of
size a most $4$. Because of the example \ref{examplskew}, such expression
can not be a polynomial. Nevertheless, it is reasonable to suggest the
following problem

\begin{prob}
Let $\left\langle n\right\rangle _{4}$ the collection of all subset $\alpha
\subseteq \left\{ 1,\ldots ,n\right\} $ with $\left\vert \alpha \right\vert $
$=2$ or $4$ and consider the polynomial ring $R=\mathbb{K}[X_{\alpha
},\alpha $ $\in \left\langle n\right\rangle _{4}]$ where $X_{\alpha }$ $%
,\alpha $ $\in \left\langle n\right\rangle _{4}$ are indeterminates. Is
there a polynomial $Q(X_{\alpha },\alpha \in \left\langle n\right\rangle
_{4})\in R$ such that
 $\prod\limits_{1\leq i<j\leq n} x_{ij}^{2}\det
\left( M\right) =Q(\det \left( M\left[ \alpha \right] \right) ,\alpha \in
\left\langle n\right\rangle _{4})$ ?
\end{prob}

Note that a positive answer to this problem combined with the Theorem \ref%
{loewy} allows to obtain our main theorem.

Finally, as the principal minors are squares of the corresponding Pfaffians,
we can strengthen the assumptions of the Problem \ref{problem2} by replacing
the principal minors by sub-pfaffians. However, two skew-symmetric matrices
with the same corresponding principal sub-pfaffians of order $2$ are equal.
So, it is interesting to consider the following problem.

\begin{prob}
Given a positive integer $n\geq 5$, what is the relationship between two
skew-symmetric matrices of orders $n$ which differ up to the sign of their
off-diagonal terms and having equal corresponding principal Pfaffian minors
of order $4$ ?
\end{prob}


\begin{thebibliography}{00}


\bibitem{BCG} A. Bouchet, W.H. Cunningham, J.F. Geelen, Principally
unimodular skew symmetric matrices, Combinatorica, 18 (1998) 461--486.

\bibitem{Bor} A. Borodin, E. Rains, Eynard--Mehta theorem, Schur process,
and their Pfaffian analogs, J. Stat. Phys. 121 (2005) 291--317.

\bibitem{BILT} A. Boussa\"{\i}ri, P. Ille, G. Lopez, S. Thomass\'{e}, The $%
C_{3}$-structure of tournaments, Discrete Math., 277 (2004), p. 29-43.

\bibitem{EHR ROZ} A. Ehrenfeucht, G. Rozenberg, Primitivity is hereditary
for 2-structures, fundamental study, Theoret. Comput. Sci. 3 (70) (1990)
343--358.

\bibitem{Engel} G. M.Engel, H.Schneider, Matrices Diagonally Similar to a
Symmetric Matrix, Linear Algebra and its Applications 29, 131-138,1980

\bibitem{Griffin1} K. Griffin and M.J. Tsatsomeros, Principal minors, part
I: a method for computing all the principal minors of a matrix, Linear
Algebra Appl., 419(1):107--124, 2006.

\bibitem{Griffin2} K. Griffin and M.J. Tsatsomeros, Principal minors, part
II: the principal minor assignment problem, Linear Algebra Appl.,
419(1):125--171, 2006.

\bibitem{HL} D.J. Hartfiel, R. Loewy, On Matrices Having equal Corresponding
Principal Minors, Linear Algebra and its Applications, 58 (1984) 147-167.

\bibitem{HOZ} O. Holtz, H. Schneider, Open problems on GKK t-matrices,
Linear Algebra Appl. 345 (2002) 263--267

\bibitem{HOLTZ07} O. Holtz, B. Sturmfels, Hyperdeterminantal relations among
symmetric principal minors, J. Algebra 316 (2007) 634--648.

\bibitem{Kenyon} R. Kenyon, R. Pemantle, Principal minors and rhombus
tilings, arXiv : 1404.1354.

\bibitem{Lin} S. Lin and B. Sturmfels, Polynomial relations among principal
minors of a $4\times 4$-matrix, J. Algebra 322 (2009), no. 11, 4121--4131.

\bibitem{Lw86} R. Loewy, Principal Minors and Diagonal Similarity of
Matrices, Linear Algebra and its Applications, 78 (1986) 23-63.

\bibitem{Oeding} L. Oeding, Set-theoretic defining equations of the variety
of principal minors of symmetric matrices, Algebra and Number Theory 5
(2011), no. 1, 75--109.

\bibitem{Rising} J. Rising, A. Kulesza, B. Taskar, An efficient algorithm
for the symmetric principal minor assignment problem, Linear Algebra and its
Applications, In Press, Corrected Proof, Available online 14 May 2014

\bibitem{Wesp} G. Wesp, A note on the spectra of certain skew-symmetric
\{1,0,-1\}-matrices, Discrete Math. 258 (2002), no. 1-3, 339-346.
\end{thebibliography}
\end{document}